\documentclass[a4paper,11pt]{amsart}
\usepackage{amsmath,amssymb,amsfonts,latexsym}

\newtheorem{theorem}{Theorem}[section]

\input{tcilatex}

\begin{document}
\title[\textbf{Extended} $q$\textbf{-Dedekind-type DC sums}]{\textbf{A note
on the} \textbf{modified} $q$\textbf{-Dedekind sums}}
\author[\textbf{S. Araci}]{\textbf{Serkan Araci}}
\address{\textbf{University of Gaziantep, Faculty of Science and Arts,
Department of Mathematics, 27310 Gaziantep, TURKEY}}
\email{\textbf{mtsrkn@hotmail.com; saraci88@yahoo.com.tr; mtsrkn@gmail.com}}
\author[\textbf{E. \c{S}en}]{\textbf{Erdo\u{g}an \c{S}en}}
\address{\textbf{Department of Mathematics, Faculty of Science and Letters,
Nam\i k Kemal University, 59030 Tekirda\u{g}, TURKEY}}
\email{\textbf{erdogan.math@gmail.com}}
\author[\textbf{M. Acikgoz}]{\textbf{Mehmet Acikgoz}}
\address{\textbf{University of Gaziantep, Faculty of Science and Arts,
Department of Mathematics, 27310 Gaziantep, TURKEY}}
\email{\textbf{acikgoz@gantep.edu.tr}}

\begin{abstract}
In the present paper, the fundamental aim is to consider a $p$-adic
continuous function for an odd prime to inside a $p$-adic $q$-analogue of
the higher order Extended Dedekind-type sums related to $q$-Genocchi
polynomials with weight $\alpha $ by using fermionic $p$-adic invariant $q$%
-integral on $%
\mathbb{Z}
_{p}$.

\vspace{2mm}\noindent \textsc{2010 Mathematics Subject Classification.}
11S80, 11B68.

\vspace{2mm}

\end{abstract}

\thanks{}
\maketitle




\section{\textbf{Introduction}}


Imagine that $p$ be a fixed odd prime number. We now start with definition
of the following notations. Let $%
\mathbb{Q}
_{p}$ be the field $p$-adic rational numbers and let $%
\mathbb{C}
_{p}$ be the completion of algebraic closure of $%
\mathbb{Q}
_{p}$.

Thus, 
\begin{equation*}
\boldsymbol{%
\mathbb{Q}
}_{p}=\left\{ x=\sum_{n=-k}^{\infty }a_{n}p^{n}:0\leq a_{n}<p\right\} \text{.%
}
\end{equation*}

Then $%
\mathbb{Z}
_{p}$ is integral domain, which is defined by 
\begin{equation*}
\boldsymbol{%
\mathbb{Z}
}_{p}=\left\{ x=\sum_{n=0}^{\infty }a_{n}p^{n}:0\leq a_{n}<p\right\}
\end{equation*}%
or%
\begin{equation*}
\boldsymbol{%
\mathbb{Z}
}_{p}=\left\{ x\in 
\mathbb{Q}
_{p}:\left\vert x\right\vert _{p}\leq 1\right\} \text{.}
\end{equation*}

We assume that $q\in 
\mathbb{C}
_{p}$ with $\left\vert 1-q\right\vert _{p}<1$ as an indeterminate. The $p$%
-adic absolute value $\left\vert .\right\vert _{p}$, is normally defined by 
\begin{equation*}
\left\vert x\right\vert _{p}=\frac{1}{p^{n}}
\end{equation*}%
where $x=p^{n}\frac{s}{t}$ with $\left( p,s\right) =\left( p,t\right)
=\left( s,t\right) =1$ and $n\in 
\mathbb{Q}
$ (for details, see [1-19]).

The $p$-adic $q$-Haar distribution is defined by Kim as follows: for any
postive integer $n$,%
\begin{equation*}
\mu _{q}\left( a+p^{n}%
\mathbb{Z}
_{p}\right) =\left( -q\right) ^{a}\frac{\left( 1+q\right) }{1+q^{p^{n}}}
\end{equation*}%
for $0\leq a<p^{n}$ and this can be extended to a measure on $%
\mathbb{Z}
_{p}$ (for details, see \cite{Kim 4}, \cite{Kim 6}, \cite{Kim 9}).

In \cite{Araci 2}, the $q$-Genocchi polynomials are defined by Araci \textit{%
et al.} as follows:%
\begin{equation}
\widetilde{G}_{n,q}^{\left( \alpha \right) }\left( x\right) =n\int_{%
\mathbb{Z}
_{p}}\left( \frac{1-q^{\alpha \left( x+\xi \right) }}{1-q^{\alpha }}\right)
^{n-1}d\mu _{q}\left( \xi \right)  \label{equation 1}
\end{equation}%
for $n\in 
\mathbb{Z}
_{+}:=\left\{ 0,1,2,3,\cdots \right\} $. We easily see that 
\begin{equation*}
\lim_{q\rightarrow 1}\widetilde{G}_{n,q}^{\left( \alpha \right) }\left(
x\right) =G_{n}\left( x\right)
\end{equation*}%
where $G_{n}\left( x\right) $ are Genocchi polynomials, which are given in
the form:%
\begin{equation*}
\sum_{n=0}^{\infty }G_{n}\left( x\right) \frac{t^{n}}{n!}=e^{tx}\frac{2t}{%
e^{t}+1},\text{ }\left\vert t\right\vert <\pi
\end{equation*}%
(for details, see \cite{Araci 2}). Taking $x=0$ into (\ref{equation 1}),
then we have $\widetilde{G}_{n,q}^{\left( \alpha \right) }\left( 0\right) :=%
\widetilde{G}_{n,q}^{\left( \alpha \right) }$ are called $q$-Genocchi
numbers with weight $\alpha $.

The $q$-Genocchi numbers and polynomials have the following identities:%
\begin{eqnarray}
\widetilde{G}_{n+1,q}^{\left( \alpha \right) } &=&\left( n+1\right) \frac{1+q%
}{\left( 1-q^{\alpha }\right) ^{n}}\sum_{l=0}^{n}\binom{n}{l}\left(
-1\right) ^{l}\frac{1}{1+q^{\alpha l+1}}\text{,}  \label{equation 2} \\
\widetilde{G}_{n+1,q}^{\left( \alpha \right) }\left( x\right) &=&\left(
n+1\right) \frac{1+q}{\left( 1-q^{\alpha }\right) ^{n}}\sum_{l=0}^{n}\binom{n%
}{l}\left( -1\right) ^{l}\frac{q^{\alpha lx}}{1+q^{\alpha l+1}}\text{,}
\label{equation 3} \\
\widetilde{G}_{n,q}^{\left( \alpha \right) }\left( x\right) &=&\sum_{l=0}^{n}%
\binom{n}{l}q^{\alpha lx}\widetilde{G}_{l,q}^{\left( \alpha \right) }\left( 
\frac{1-q^{\alpha x}}{1-q^{\alpha }}\right) ^{n-l}\text{.}
\label{equation 4}
\end{eqnarray}

Additionally, for $d$ odd natural number, we have 
\begin{equation}
\widetilde{G}_{n,q}^{\left( \alpha \right) }\left( dx\right) =\left( \frac{%
1+q}{1+q^{d}}\right) \left( \frac{1-q^{\alpha d}}{1-q^{\alpha }}\right)
^{n-1}\sum_{a=0}^{d-1}q^{a}\left( -1\right) ^{a}\widetilde{G}_{n,q}^{\left(
\alpha \right) }\left( x+\frac{a}{d}\right) \text{,}  \label{equation 5}
\end{equation}

(for details about this subject, see \cite{Araci 2}).

For any positive integer $h,k$ and $m$, Dedekind-type DC sums are given by
Kim in \cite{Kim 1}, \cite{Kim 2} and \cite{Kim 3} as follows:%
\begin{equation*}
S_{m}\left( h,k\right) =\sum_{M=1}^{k-1}\left( -1\right) ^{M-1}\frac{M}{k}%
\overline{E}_{m}\left( \frac{hM}{k}\right)
\end{equation*}%
where $\overline{E}_{m}\left( x\right) $ are the $m$-th periodic Euler
function.

In 2011, Taekyun Kim added a weight to $q$-Bernoulli polynomials in \cite%
{Kim 8}. He derived not only new but also ineteresting properties for
weighted $q$-Bernoulli polynomials. After, many mathematicians, by utilizing
from Kim's paper \cite{Kim 8}, have introduced a new concept in Analytic
numbers theory as weighted $q$-Bernoulli, weighted $q$-Euler, weighted $q$%
-Genocchi polynomials in \cite{Kim 9}, \cite{Araci 1}, \cite{Araci 2}, \cite%
{Araci 3}, \cite{Araci 6} and \cite{Araci 5}. Also, the generating function
of weighted $q$-Genocchi polynomials was introduced by Araci \textit{et al}.
in \cite{Araci 2}. They also derived several arithmetic properties for
weighted $q$-Genocchi polynomials.

Kim has given some interesting properties for Dedekind-type DC sums. He
firstly considered a $p$-adic continuous function for an odd prime number to
contain a $p$-adic $q$-analogue of the higher order Dedekind-type DC sums $%
k^{m}S_{m+1}\left( h,k\right) $ in \cite{Kim 2}.

By the same motivation, we, by using $p$-adic invariant $q$-integral on $%
\mathbb{Z}
_{p}$, shall get weighted $p$-adic $q$-analogue of the higher order
Dedekind-type DC sums $k^{m}S_{m+1}\left( h,k\right) $.

\section{\textbf{Extended }$q$\textbf{-Dedekind-type Sums in connection with 
}$q$\textbf{-Genocchi polynomials with weight }$\protect\alpha $}

If $x$ is a $p$-adic integer, then $w\left( x\right) $ is the unique
solution of $w\left( x\right) =w\left( x\right) ^{p}$ that is congruent to $x%
\func{mod}p$. It can also be defined by%
\begin{equation*}
w\left( x\right) =\lim_{n\rightarrow \infty }x^{p^{n}}\text{.}
\end{equation*}

The multiplicative group of $p$-adic units is a product of the finite group
of roots of unity, and a group isomorphic to the $p$-adic integers. The
finite group is cylic of order $p-1$ or $2$, as $p$ is odd or even,
respectively, and so it is isomorphic. Actually, the teichm\"{u}ller
character gives a canonical isomorphism between these two groups.

Let $w$ be the $Teichm\ddot{u}ller$ character ($\func{mod}p$). For $x\in 
\mathbb{Z}
_{p}^{\ast }$ $:=%
\mathbb{Z}
_{p}/p%
\mathbb{Z}
_{p}$, set%
\begin{equation*}
\left\langle x:q\right\rangle =w^{-1}\left( x\right) \left( \frac{1-q^{x}}{%
1-q}\right) \text{.}
\end{equation*}

Let $a$ and $N$ be positive integers with $\left( p,a\right) =1$ and $p\mid
N $. We now introduce the following%
\begin{equation*}
\widetilde{E}_{q}^{\left( \alpha \right) }\left( s,a,N:q^{N}\right)
=w^{-1}\left( a\right) \left\langle x:q^{\alpha }\right\rangle
^{s}\sum_{j=0}^{\infty }\binom{s}{j}q^{\alpha aj}\left( \frac{1-q^{\alpha N}%
}{1-q^{\alpha a}}\right) ^{j}\widetilde{G}_{j,q^{N}}^{\left( \alpha \right) }%
\text{.}
\end{equation*}

In particular, if $m+1\equiv 0(\func{mod}p-1)$, then we have%
\begin{eqnarray*}
\widetilde{E}_{q}^{\left( \alpha \right) }\left( m,a,N:q^{N}\right)
&=&\left( \frac{1-q^{\alpha a}}{1-q^{\alpha }}\right) ^{m}\sum_{j=0}^{m}%
\binom{m}{j}q^{\alpha aj}\widetilde{G}_{j,q^{N}}^{\left( \alpha \right)
}\left( \frac{1-q^{\alpha N}}{1-q^{\alpha a}}\right) ^{j} \\
&=&\left( \frac{1-q^{\alpha N}}{1-q^{\alpha }}\right) ^{m}\int_{%
\mathbb{Z}
_{p}}\left( \frac{1-q^{\alpha N\left( \xi +\frac{a}{N}\right) }}{1-q^{\alpha
N}}\right) ^{m}d\mu _{q^{N}}\left( \xi \right) \text{.}
\end{eqnarray*}

Then, $\widetilde{E}_{q}^{\left( \alpha \right) }\left( m,a,N:q^{N}\right) $
is a continuous $p$-adic extension of 
\begin{equation*}
\left( \frac{1-q^{\alpha N}}{1-q^{\alpha }}\right) ^{m}\frac{\widetilde{G}%
_{m+1,q^{N}}^{\left( \alpha \right) }\left( \frac{a}{N}\right) }{m+1}\text{.}
\end{equation*}

Suppose that $\left[ .\right] $ be the Gauss' symbol and let $\left\{
x\right\} =x-\left[ x\right] $. Thus, we are now ready to treat $q$%
-extension of the higher order Dedekind-type DC sums $\widetilde{S}%
_{m,q}^{\left( \alpha \right) }\left( h,k:q^{l}\right) $ in the form: 
\begin{equation*}
\widetilde{S}_{m,q}^{\left( \alpha \right) }\left( h,k:q^{l}\right)
=\sum_{M=1}^{k-1}\left( -1\right) ^{M-1}\left( \frac{1-q^{\alpha M}}{%
1-q^{\alpha k}}\right) \int_{%
\mathbb{Z}
_{p}}\left( \frac{1-q^{\alpha \left( l\xi +l\left\{ \frac{hM}{k}\right\}
\right) }}{1-q^{\alpha l}}\right) ^{m}d\mu _{q^{l}}\left( \xi \right) \text{.%
}
\end{equation*}

If $m+1\equiv 0\left( \func{mod}p-1\right) $%
\begin{eqnarray*}
&&\left( \frac{1-q^{\alpha k}}{1-q^{\alpha }}\right)
^{m+1}\sum_{M=1}^{k-1}\left( -1\right) ^{M-1}\left( \frac{1-q^{\alpha M}}{%
1-q^{\alpha k}}\right) \int_{%
\mathbb{Z}
_{p}}\left( \frac{1-q^{\alpha k\left( \xi +\frac{hM}{k}\right) }}{%
1-q^{\alpha k}}\right) ^{m}d\mu _{q^{k}}\left( \xi \right) \\
&=&\sum_{M=1}^{k-1}\left( -1\right) ^{M-1}\left( \frac{1-q^{\alpha M}}{%
1-q^{\alpha }}\right) \left( \frac{1-q^{\alpha k}}{1-q^{\alpha }}\right)
^{m}\int_{%
\mathbb{Z}
_{p}}\left( \frac{1-q^{\alpha k\left( \xi +\frac{hM}{k}\right) }}{%
1-q^{\alpha k}}\right) ^{m}d\mu _{q^{k}}\left( \xi \right)
\end{eqnarray*}%
where $p\mid k$, $\left( hM,p\right) =1$ for each $M$. Via the equation (\ref%
{equation 1}), we easily state the following 
\begin{gather}
\left( \frac{1-q^{\alpha k}}{1-q^{\alpha }}\right) ^{m+1}\widetilde{S}%
_{m,q}^{\left( \alpha \right) }\left( h,k:q^{k}\right)  \label{equation 6} \\
=\sum_{M=1}^{k-1}\left( \frac{1-q^{\alpha M}}{1-q^{\alpha }}\right) \left( 
\frac{1-q^{\alpha k}}{1-q^{\alpha }}\right) ^{m}\left( -1\right) ^{M-1}\int_{%
\mathbb{Z}
_{p}}\left( \frac{1-q^{\alpha k\left( \xi +\frac{hM}{k}\right) }}{%
1-q^{\alpha k}}\right) ^{m}d\mu _{q^{k}}\left( \xi \right)  \notag \\
=\sum_{M=1}^{k-1}\left( -1\right) ^{M-1}\left( \frac{1-q^{\alpha M}}{%
1-q^{\alpha }}\right) \widetilde{E}_{q}^{\left( \alpha \right) }\left(
m,\left( hM\right) _{k}:q^{k}\right)  \notag
\end{gather}%
where $(hM)_{k}$ denotes the integer $x$ such that $0\leq x<n$ and $x\equiv
\alpha \left( \func{mod}k\right) $.

It is simple to show the following: 
\begin{gather}
\int_{%
\mathbb{Z}
_{p}}\left( \frac{1-q^{\alpha \left( x+\xi \right) }}{1-q^{\alpha }}\right)
^{k}d\mu _{q}\left( \xi \right)  \label{equation 7} \\
=\left( \frac{1-q^{\alpha m}}{1-q^{\alpha }}\right) ^{k}\frac{1+q}{1+q^{m}}%
\sum_{i=0}^{m-1}\left( -1\right) ^{i}\int_{%
\mathbb{Z}
_{p}}\left( \frac{1-q^{\alpha m\left( \xi +\frac{x+i}{m}\right) }}{%
1-q^{\alpha m}}\right) ^{k}d\mu _{q^{m}}\left( \xi \right) \text{.}  \notag
\end{gather}

Due to (\ref{equation 6}) and (\ref{equation 7}), we easily obtain%
\begin{gather}
\left( \frac{1-q^{\alpha N}}{1-q^{\alpha }}\right) ^{m}\int_{%
\mathbb{Z}
_{p}}\left( \frac{1-q^{\alpha N\left( \xi +\frac{a}{N}\right) }}{1-q^{\alpha
N}}\right) ^{m}d\mu _{q^{N}}\left( \xi \right)  \label{equation 8} \\
=\frac{1+q^{N}}{1+q^{Np}}\sum_{i=0}^{p-1}\left( -1\right) ^{i}\left( \frac{%
1-q^{\alpha Np}}{1-q^{\alpha }}\right) ^{m}\int_{%
\mathbb{Z}
_{p}}\left( \frac{1-q^{\alpha pN\left( \xi +\frac{a+iN}{pN}\right) }}{%
1-q^{\alpha pN}}\right) ^{m}d\mu _{q^{pN}}\left( \xi \right) \text{.}  \notag
\end{gather}

Thanks to (\ref{equation 6}), (\ref{equation 7}) and (\ref{equation 8}), we
discover the following $p$-adic integration: 
\begin{equation*}
\widetilde{E}_{q}^{\left( \alpha \right) }\left( s,a,N:q^{N}\right) =\frac{%
1+q^{N}}{1+q^{Np}}\sum_{\underset{a+iN\neq 0(\func{mod}p)}{0\leq i\leq p-1}%
}\left( -1\right) ^{i}\widetilde{E}_{q}^{\left( \alpha \right) }\left(
s,\left( a+iN\right) _{pN},p^{N}:q^{pN}\right) \text{.}
\end{equation*}

On the other hand,%
\begin{gather*}
\widetilde{E}_{q}^{\left( \alpha \right) }\left( m,a,N:q^{N}\right) =\left( 
\frac{1-q^{\alpha N}}{1-q^{\alpha }}\right) ^{m}\int_{%
\mathbb{Z}
_{p}}\left( \frac{1-q^{\alpha N\left( \xi +\frac{a}{N}\right) }}{1-q^{\alpha
N}}\right) ^{m}d\mu _{q^{N}}\left( \xi \right) \\
-\left( \frac{1-q^{\alpha Np}}{1-q^{\alpha }}\right) ^{m}\int_{%
\mathbb{Z}
_{p}}\left( \frac{1-q^{\alpha pN\left( \xi +\frac{a+iN}{pN}\right) }}{%
1-q^{\alpha pN}}\right) ^{m}d\mu _{q^{pN}}\left( \xi \right)
\end{gather*}%
where $\left( p^{-1}a\right) _{N}$ denotes the integer $x$ with $0\leq x<N$, 
$px\equiv a\left( \func{mod}N\right) $ and $m$ is integer with $m+1\equiv 0(%
\func{mod}p-1)$. Therefore, we can state the following%
\begin{gather*}
\sum_{M=1}^{k-1}\left( -1\right) ^{M-1}\left( \frac{1-q^{\alpha M}}{%
1-q^{\alpha }}\right) \widetilde{E}_{q}^{\left( \alpha \right) }\left(
m,hM,k:q^{k}\right) \\
=\left( \frac{1-q^{\alpha k}}{1-q^{\alpha }}\right) ^{m+1}\widetilde{S}%
_{m,q}^{\left( \alpha \right) }\left( h,k:q^{k}\right) -\left( \frac{%
1-q^{\alpha k}}{1-q^{\alpha }}\right) ^{m+1}\left( \frac{1-q^{\alpha kp}}{%
1-q^{\alpha k}}\right) \widetilde{S}_{m,q}^{\left( \alpha \right) }\left(
\left( p^{-1}h\right) ,k:q^{pk}\right)
\end{gather*}%
where $p\nmid k$ and $p\nmid hm$ for each $M$. Thus, we obtain the following
definition, which seems interesting for further studying in theory of
Dedekind sums.

\begin{definition}
Let $h,k$ be positive integer with $\left( h,k\right) =1$, $p\nmid k$. For $%
s\in 
\mathbb{Z}
_{p},$ we define $p$-adic Dedekind-type DC sums as follows:  
\begin{equation*}
\widetilde{S}_{p,q}^{\left( \alpha \right) }\left( s:h,k:q^{k}\right)
=\sum_{M=1}^{k-1}\left( -1\right) ^{M-1}\left( \frac{1-q^{\alpha M}}{%
1-q^{\alpha }}\right) \widetilde{E}_{q}^{\left( \alpha \right) }\left(
m,hM,k:q^{k}\right) \text{.}
\end{equation*}
\end{definition}

As a result of the above definition, we derive the following theorem.

\begin{theorem}
For $m+1\equiv 0(\func{mod}p-1)$ and $\left( p^{-1}a\right) _{N}$ denotes
the integer $x$ with $0\leq x<N$, $px\equiv a\left( \func{mod}N\right) $,
then, we have 
\begin{gather*}
\widetilde{S}_{p,q}^{\left( \alpha \right) }\left( s:h,k:q^{k}\right)
=\left( \frac{1-q^{\alpha k}}{1-q^{\alpha }}\right) ^{m+1}\widetilde{S}%
_{m,q}^{\left( \alpha \right) }\left( h,k:q^{k}\right) \\
-\left( \frac{1-q^{\alpha k}}{1-q^{\alpha }}\right) ^{m+1}\left( \frac{%
1-q^{\alpha kp}}{1-q^{\alpha k}}\right) \widetilde{S}_{m,q}^{\left( \alpha
\right) }\left( \left( p^{-1}h\right) ,k:q^{pk}\right) \text{.}
\end{gather*}
\end{theorem}

In the special case $\alpha =1$, our applications in theory of Dedekind sums
resemble Kim's results in \cite{Kim 2}. These results seem to be interesting
for further studies in \cite{Kim 1}, \cite{Kim 3} and \cite{Simsek}.



\end{document}